\documentclass[english]{svjour3}
\usepackage[latin9]{inputenc}
\usepackage{geometry}
\usepackage{babel}
\usepackage{url}
\usepackage{amsmath}
\usepackage{graphicx}
\usepackage{esint}
\usepackage[numbers]{natbib}
\usepackage[unicode=true,pdfusetitle,
 bookmarks=true,bookmarksnumbered=true,bookmarksopen=false,
 breaklinks=false,pdfborder={0 0 0},backref=false,colorlinks=false]
 {hyperref}
\usepackage{breakurl}

\makeatletter

\providecommand{\tabularnewline}{\\}
\newcommand{\lyxdot}{.}

\numberwithin{equation}{section}
\numberwithin{figure}{section}
  \newenvironment{svmultproof}{\begin{proof}}{\qed\end{proof}}

\usepackage{babel}

\usepackage{babel}

\usepackage{amsfonts}
\usepackage{graphicx}
 \usepackage{layout}
\usepackage{caption}

\makeatother

\begin{document}

\title{Bernstein modal basis: application to the spectral Petrov-Galerkin
method for fractional partial differential equations}

\titlerunning{Bernstein modal basis}

\author{M. Jani \and E. Babolian \and S. Javadi}

\institute{M. Jani \and E. Babolian \and S. Javadi \at Department of Mathematics,
Faculty of Mathematical Sciences and Computer, Kharazmi University,
Tehran, Iran\\
\email{mostafa.jani@gmail.com babolian@khu.ac.ir javadi@khu.ac.ir.}
\\
}
\maketitle
\begin{abstract}
In the spectral Petrov-Galerkin methods, the trial and test functions
are required to satisfy particular boundary conditions. By a suitable
linear combination of orthogonal polynomials, a basis, that is called
the modal basis, is obtained. In this paper, we extend this idea to
the non-orthogonal dual Bernstein polynomials. A compact general formula
is derived for the modal basis functions based on dual Bernstein polynomials.
Then, we present a Bernstein-spectral Petrov-Galerkin method for a
class of time fractional partial differential equations with Caputo
derivative. It is shown that the method leads to banded sparse linear
systems for problems with constant coefficients. Some numerical examples
are provided to show the efficiency and the spectral accuracy of the
method.

\keywords{Bernstein polynomials \and Petrov-Galerkin \and Dual Bernstein polynomials
\and Fractional partial differential equations \and Modal basis} \subclass{35R11 \and 65M70\and 41A10 \and 65M22 \and 76M22}
\end{abstract}

\section{Introduction}

Due to the interesting features like shape preserving \citep{carnicer},
optimal stability \citep{farouki2-1996}, etc., Bernstein polynomials
are commonly used in computer aided geometric design (CAGD) for approximating
curves and surfaces and designing computer fonts \citep{farin2002}.
They have been applied in popular programs such as Adobe's Illustrator,
Flash and Postscript in the form of Bézier curves \citep{Lawson2012}.

Bernstein polynomials have also been implemented for solving differential,
integro-differential and fractional differential equations \citep{Behiry,DohaInt,javadi,saadatmandi}.
However, they are not orthogonal, leading to dense linear systems.
The dual Bernstein polynomials (DBP) were explicitly presented by
Juttler in 1998 \citep{juttler}. To the best of our knowledge, they
have been discussed only from CAGD point of view \citep{Lewanowicz,wozny}.

The main purpose of this work is to derive a new polynomial basis
by using the DBPs that can be used with the Petrov-Galerkin formulation
for the boundary value problems of any order. It can also be used
for solving fractional differential equations. We present a Bernstein-spectral
Petrov-Galerkin method for a class of time fractional partial differential
equations. It is shown that the method leads to banded matrices for
problems with constant coefficients, saving in the computational costs
for a desired accuracy. The spectral methods utilize high order basis
functions, typically the orthogonal polynomials that are the solutions
of the Sturm-Liouville equation and they are known to have spectral
accuracy, i.e., having convergence speed faster than the methods with
fixed polynomial rate of convergence like the finite element and finite
difference methods, for problems with smooth solution. 

Fractional PDEs play a key role in modeling some physical phenomena
such as particle transport process in anomalous diffusion which has
applications in semiconductors, finance, electrochemistry, etc. \citep{Bhrawy,Gao,Goychuk,Izadkhah}.
The Caputo temporal fractional derivative of $u\left(x,t\right)$
is defined as 
\begin{eqnarray}
 &  & {\partial{}_{t}^{\alpha}}u\left(x,t\right)=\frac{1}{\Gamma\left(1-\alpha\right)}\int_{0}^{t}{\frac{1}{\left(t-s\right){}^{\alpha}}\frac{\partial u\left(x,s\right)}{\partial s}ds},\quad0<\alpha<1.\label{eq:FracDiff}
\end{eqnarray}

The paper is organized as follows. Section \ref{sec:Bern=000026Dual}
gives some preliminaries of Bernstein and DBPs. Some new aspects of
these polynomials and an interesting formula for the modal basis functions
are presented in Section \ref{sec:Modal-basis-functions}. In Section
\ref{sec:dual-Petrov-Gale=000026Error}, a Bernstein-spectral Petrov-Galerkin
method is developed for a class of time fractional differential equations.
The efficiency and spectral accuracy of the method are illustrated
via numerical examples in Section \ref{sec:NumericalEx}.

\section{\label{sec:Bern=000026Dual}Bernstein polynomials and DBPs}

Bernstein basis polynomials of degree $N$ over the unit interval
$I=[0,1]$ are defined by
\begin{eqnarray}
\phi_{i}:=B_{i,N}(x) & = & \binom{N}{i}x^{i}(1-x)^{N-i},\quad0\leq i\leq N.\label{eq:Bpolys}
\end{eqnarray}
We adopt the convention $\phi_{i}(x)\equiv0$ for $i<0$ and $i>N$.
The set $\{\phi_{i}(x):0\leq i\leq N\}$ forms a basis for $P_{N}$,
the space of polynomials with degree not exceeding $N$. It enjoys
interesting properties facilitating the numerical implementation.
They possess the end-point interpolation property  \citep{janiNA}
\begin{eqnarray}
 &  & \phi_{i}(0)=\delta_{i,0},\quad\phi_{i}(1)=\delta_{i,N},\label{eq:End-pointDelta}\\
 &  & \phi_{i}^{(p)}(0)=\frac{(-1)^{i+p}N!}{(N-p)!}\binom{p}{i},\quad p\leq N,\label{eq:Value of higher orders at 0 and 1}\\
 &  & \phi_{i}^{(p)}(1)=\frac{(-1)^{N-i}N!}{(N-p)!}\binom{p}{N-i},\quad p\leq N,
\end{eqnarray}
for $0\leq i\leq N.$ Especially, for $n\in\mathbb{N}$, $n\le N$,
the polynomials $\phi_{i}(x)$, $\lfloor\frac{n}{2}\rfloor\leq i\leq N-\lfloor\frac{n+1}{2}\rfloor$,
satisfy the end-point conditions 
\begin{eqnarray}
 &  & \phi_{i}^{(p)}(0)=\phi_{i}^{(p)}(1)=0,\quad0\leq p\leq\lfloor\frac{n}{2}\rfloor-1,\label{eq:CondPhi}
\end{eqnarray}
where $\left\lfloor \cdot\right\rfloor $ indicates the floor function.
Moreover, when $n$ is an odd integer, we get 
\begin{eqnarray}
 &  & \phi_{i}^{(\lfloor\frac{n}{2}\rfloor)}(1)=0,\quad\lfloor\frac{n}{2}\rfloor\leq i\leq N-\lfloor\frac{n+1}{2}\rfloor.\label{eq:CondPhii}
\end{eqnarray}
We will use (\ref{eq:CondPhi}) and (\ref{eq:CondPhii}) to introduce
a basis for the solution fractional partial differential equations. 

\begin{flushleft}
The derivative of Bernstein polynomials satisfies a three-term recurrence
formula \citep{janiNA} 
\begin{eqnarray}
 &  & \phi_{i}^{\prime}\left(x\right)=\left(N-i+1\right)\phi_{i-1}\left(x\right)-\left(N-2i\right)\phi_{i}\left(x\right)-\left(i+1\right)\phi_{i+1}\left(x\right),\quad0\leq i\leq N.\label{eq:ThreeTermDer}
\end{eqnarray}
The dual Bernstein polynomials given by 
\begin{eqnarray}
 & \tilde{\psi}_{i} & (x)=\sum_{j=0}^{N}c_{ij}\phi_{j}(x),\label{eq:dualPolys}
\end{eqnarray}
with the coefficients 
\begin{eqnarray*}
 &  & c_{ij}=\frac{(-1)^{i+j}}{\binom{N}{i}\binom{N}{j}}\sum_{r=0}^{\min(i,j)}(2r+1)\binom{N+r+1}{N-i}\binom{N-r}{N-i}\binom{N+r+1}{N-j}\binom{N-r}{N-j},
\end{eqnarray*}
provide the following biorthogonality system 
\begin{eqnarray}
 &  & (\phi_{i},\tilde{\psi}_{j})=\delta_{i,j},\quad0\leq i,j\leq N.\label{eq:Biorthog}
\end{eqnarray}
with the standard $L^{2}$ inner product $(f,g)=\int_{I}{fg}dx$ \citep{juttler}.
The matrix $\mathbf{C}=[c_{i,j}:0\leq i,j\leq N]$ is bisymmetric,
i.e., $c_{i,j}=c_{j,i}=c_{N-i,N-j}$. It is also seen that 
\begin{eqnarray}
 &  & \sum_{i=0}^{N}c_{i,j}=\sum_{j=0}^{N}c_{i,j}=N+1,\quad0\leq i,j\leq N.\label{eq:SumOfCentries}
\end{eqnarray}
This can be proved by (\ref{eq:dualPolys}) and changing the order
of the double summation.
\par\end{flushleft}

\section{\label{sec:Modal-basis-functions}Modal basis functions}

Using a suitable linear combination of a known orthogonal basis, typically
the Jacobi polynomial basis, one may form a basis for the spectral
Petrov-Galerkin method (see e.g. \citep{GoubetShenDual,YuanShenDual}).
We extend this idea to the non-orthogonal dual Bernstein polynomials.
In this section, a compact formula for the modal basis functions is
presented for an arbitrary order boundary value problem (BVP).

Consider a BVP of order $n$ with boundary conditions 
\begin{eqnarray}
 &  & \begin{array}{cc}
v^{(i)}(0)=v^{(i)}(1)=0,\quad0\leq i\leq\frac{n-2}{2},\hfill\quad\quad\,\, & \mathrm{for\medspace}n\medspace\mathrm{even},\\
v^{(i)}(0)=v^{(i)}(1)=0,\quad0\leq i\leq\frac{n-3}{2},\medspace\mathrm{and}\medspace v^{(\frac{n-1}{2})}(1)=0, & \mathrm{for}\medspace n\medspace\mathrm{odd}.
\end{array}\label{eq:ConditionsOfTrial}
\end{eqnarray}
There is no loss of generality in assuming homogeneous conditions.
We define the \textit{trial space} from which we seek an approximate
solution of the problem as $V_{N}^{0,n}=\{v\in P_{N}:v\medspace\mathrm{satisfies}\medspace\mathrm{the}\medspace\mathrm{conditions}\medspace(\ref{eq:ConditionsOfTrial})\}.$
We also define the \textit{test space} $W_{N}^{0,n}$, as the set
of polynomials $w$ in $P_{N}$ such that 
\begin{eqnarray}
 &  & \begin{array}{cc}
w^{(i)}(0)=w^{(i)}(1)=0,\quad0\leq i\leq\frac{n-2}{2},\hfill\quad\quad\,\, & \mathrm{for\medspace}n\medspace\mathrm{even},\\
w^{(i)}(0)=w^{(i)}(1)=0,\quad0\leq i\leq\frac{n-3}{2},\medspace\mathrm{and}\medspace w^{(\frac{n-1}{2})}(0)=0, & \mathrm{for}\medspace n\medspace\mathrm{odd}.
\end{array}\label{eq:ConditionsOfTest}
\end{eqnarray}
A basis in $W_{N}^{0,n}$ is chosen to serve as the \textit{test functions}
in Petrov-Galerkin formulation of the problem. Note that $\mathrm{dim}V_{N}^{0,n}=\mathrm{dim}W_{N}^{0,n}=N-n+1.$
Also, $V_{N}^{0,n}=W_{N}^{0,n}$ when $n$ is even. 

From (\ref{eq:CondPhi})-(\ref{eq:CondPhii}), it is seen that the
set 
\begin{eqnarray}
 &  & \{\phi_{i}(x):\left\lfloor \frac{n}{2}\right\rfloor \leq i\leq N-\lfloor\frac{n+1}{2}\rfloor\},\label{eq:trialFunctions}
\end{eqnarray}
forms a basis for $V_{N}^{0,n}$. Before presenting a basis for $W_{N}^{0,n}$,
we provide the following results for DBPs.
\begin{lemma}
\citep{janiCAMWA} Set $\alpha_{i,0}:=-(-1)^{i}(N+1)\binom{N+1}{i+1}+N\delta_{i,0}+\delta_{i,1}$
for $0\leq i\leq N.$ Then,
\begin{eqnarray}
\tilde{\psi}_{i}^{\prime}(x) & = & \alpha_{i,0}\tilde{\psi}_{0}\left(x\right)+(1-\delta_{i,1})i\tilde{\psi}_{i-1}\left(x\right)+(1-\delta_{i,0})(1-\delta_{i,N})\left(N-2i\right)\tilde{\psi}_{i}\left(x\right)\label{eq:der of DBPs}\\
 &  & \hfill-(1-\delta_{i,N-1})\left(N-i\right)\tilde{\psi}_{i+1}\left(x\right)-\alpha_{N-i,0}\tilde{\psi}_{N}\left(x\right),\nonumber 
\end{eqnarray}
where we set $\tilde{\psi}_{i}\equiv0$ for $i<0$ and $i>N.$\end{lemma}
\begin{proposition}
\label{Prop:dualproperties}The following statements hold for $0\leq i\leq N$
and $x\in I=[0,1]$:

\begin{minipage}[t]{1\columnwidth}%
\begin{tabular}{ll}
(a) $\tilde{\psi}_{N-i}\left(x\right)=\tilde{\psi}_{i}\left(1-x\right),$ & (b) $\sum_{j=0}^{N}{\tilde{\psi}_{j}\left(x\right)}=N+1,$\tabularnewline
(c) $\int_{0}^{1}\tilde{\psi}_{i}\left(x\right)dx=1,$ & (d) $\tilde{\psi}_{i}^{(p)}(0)=\frac{(-1)^{p}N!}{(N-p)!}\sum_{r=0}^{p}(-1)^{r}c_{i,r}\binom{p}{r}.$\tabularnewline
\end{tabular}%
\end{minipage}\end{proposition}
\begin{svmultproof}
The first statement follows from definition (\ref{eq:dualPolys})
and the similar relation $\phi_{N-i}(x)=\phi_{i}(1-x)$. From (\ref{eq:dualPolys})
and (\ref{eq:SumOfCentries}) we get the following that proves (b):
\begin{eqnarray*}
\sum_{i=0}^{N}{\tilde{\psi}_{i}\left(x\right)} & = & \sum_{i=0}^{N}{\sum_{j=0}^{N}c_{i,j}\phi_{j}\left(x\right)}=\sum_{j=0}^{N}{\phi_{j}\left(x\right)\sum_{i=0}^{N}c_{i,j}}=N+1.
\end{eqnarray*}
The statement (c) follows from the fact $\int_{0}^{1}\phi_{i}(x)dx=\frac{1}{N+1}$
and (\ref{eq:SumOfCentries}). (d) is derived by (\ref{eq:Value of higher orders at 0 and 1}).
\end{svmultproof}

The following theorem gives a formula for the modal basis functions
for the test space $W_{N}^{0,n}$.
\begin{theorem}
\label{TheoremModal}Let $n<N$. The following polynomials form a
basis for $W_{N}^{0,n}$. 
\begin{eqnarray}
 &  & \psi_{i}\left(x\right)=\sum_{j=0}^{n}{a_{i,j}^{n}\tilde{\psi}_{i+j}\left(x\right)},\quad0\leq i\leq N-n,\label{eq:Modal}\\
 &  & a_{i,j}^{n}=\frac{\binom{n}{j}(i+j+[\frac{n+1}{2}])!(N-i-j+[\frac{n}{2}])!}{(i+[\frac{n+1}{2}])!(N-i+[\frac{n}{2}])!},\quad0\leq j\leq n.
\end{eqnarray}
\end{theorem}
\begin{svmultproof}
The leading coefficient in (\ref{eq:Modal}) is $a_{i,0}^{n}=1$,
$\tilde{\psi_{i}}$'s are linearly independent and the number of $\psi_{i}$'s
is equal to$\dim W_{N}^{0,n}=N-n+1$. It is thus sufficient to prove
that $\psi_{i}\in W_{N}^{0,n}$ for $0\leq i\leq N-n$. If $n$ is
even, it is 
\begin{eqnarray*}
 &  & \psi_{i}^{(p)}(0)=\psi_{i}^{(p)}(1)=0,\quad0\leq p\leq\frac{n-2}{2}.
\end{eqnarray*}
To do this, using Proposition \ref{Prop:dualproperties}, we have
\begin{eqnarray*}
\psi_{i}^{(p)}(0) & = & \sum_{j=0}^{n}{a_{i,j}^{n}\tilde{\psi}_{i+j}^{(p)}\left(0\right)}\\
 & = & \frac{(-1)^{p}N!}{(N-p)!}\sum_{j=0}^{n}{a_{i,j}^{n}\sum_{r=0}^{p}(-1)^{r}c_{i+j,r}\binom{p}{r}}\\
 & = & \frac{(-1)^{p}N!}{(N-p)!}\sum_{r=0}^{p}{(-1)^{r}\binom{p}{r}\sum_{j=0}^{n}a_{i,j}^{n}c_{i+j,r}}.
\end{eqnarray*}
With some manipulations, it is seen that the inner summation vanishes
for $0\leq r\leq p$, hence the proof is completed. The proof for
odd $n$ is done similarly.
\end{svmultproof}

For example, the modal basis functions (\ref{eq:Modal}) for $W_{N}^{0,2},W_{N}^{0,3}$
and $W_{N}^{0,4}$ are written as 
\begin{flalign}
\psi_{i} & =\tilde{\psi}_{i}+\frac{i+2}{N-i+1}(2\tilde{\psi}_{i+1}+\frac{i+3}{N-i}\tilde{\psi}_{i+2}),\quad0\leq i\leq N-2,\label{eq:base for n=00003D2}\\
\psi_{i} & =\tilde{\psi}_{i}+\frac{i+2}{N-i+2}(3\tilde{\psi}_{i+1}+\frac{i+3}{N-i+1}(3\tilde{\psi}_{i+2}+\frac{i+4}{N-i}\tilde{\psi}_{i+3})),\quad0\leq i\leq N-3,\label{eq:base for n=00003D3}\\
\psi_{i} & =\tilde{\psi}_{i}+\frac{i+3}{N-i+2}(4\tilde{\psi}_{i+1}+\frac{i+4}{N-i+1}(6\tilde{\psi}_{i+2}+\frac{i+5}{N-i}(4\tilde{\psi}_{i+3}+\frac{i+6}{N-i-1}\tilde{\psi}_{i+4}))),0\leq i\leq N-4,\label{eq:base for n=00003D4}
\end{flalign}
respectively. These are used for second, third and fourth order differential
equations with conditions (\ref{eq:ConditionsOfTrial}), respectively.
As in the finite element method, the advantage of using such a basis
utilizing the neighboring functions lies in the fact that it minimizes
the interactions of basis functions in frequency space \citep{Shen}.

\section{\label{sec:dual-Petrov-Gale=000026Error} The Bernstein-spectral
Petrov-Galerkin method}

In this section, a Petrov-Galerkin method based on the modal basis
functions introduced in Theorem \ref{TheoremModal} is presented for
the time-fractional differential equation
\begin{eqnarray}
 &  & {\partial{}_{t}^{\alpha}}u(x,t)=\sum_{r=0}^{n}{b_{r}(x,t)\partial_{x}^{r}u(x,t)}+s(x,t),\quad(x,t)\in\Omega\times(0,T],\label{eq:main}
\end{eqnarray}
with $\Omega=(0,1)$, $0<\alpha\leq1$, the source term $s$, the
initial condition $u(x,0)=g(x)$ and $n$ boundary conditions (\ref{eq:ConditionsOfTrial})
in which $v(\cdot):=u(\cdot,t)$. $b_{i}$'s are given functions and
${\partial{}_{t}^{\alpha}}u$ is the Caputo derivative defined by
(\ref{eq:FracDiff}). Equation (\ref{eq:main}) includes some important
problems in science and engineering like the fractional advection-dispersion,
the anomalous diffusion, etc. \citep{Goychuk}

Let $\tau=\frac{T}{M}$ be the time step length, $t_{k}=k\tau$ and
$u^{k}\left(x\right):=u\left(x,t_{k}\right),$ $0\leq k\leq M$. The
Caputo derivative may be discretized at $t=t_{k+1},\,k\geq0$ by the
well-known L1 approximation \citep{Gao} 
\begin{eqnarray}
{\partial{}_{t}^{\alpha}}u\left(x,t_{k+1}\right)=\mu_{\tau}^{\alpha}\sum_{j=0}^{k}a_{k,j}^{\alpha}\left(u\left(x,t_{j+1}\right)-u\left(x,t_{j}\right)\right)+r_{\tau}^{k+1},\label{eq:L1Approx}
\end{eqnarray}
where $\mu_{\tau}^{\alpha}=\frac{1}{\tau^{\alpha}\Gamma\left(2-\alpha\right)}$,
$a_{kj}^{\alpha}=\left(k+1-j\right)^{1-\alpha}-\left(k-j\right)^{1-\alpha}$
and $\left|r_{\tau}^{k+1}\right|\leq\tilde{c}_{u}\tau^{2-\alpha}$
in which $\tilde{c}_{u}$ depends only on $u$ \citep{Gao}\@. Using
(\ref{eq:L1Approx}) in (\ref{eq:main}), we get 
\begin{eqnarray}
 &  & \mu_{\tau}^{\alpha}u^{k+1}\left(x\right)-\sum_{r=0}^{n}{b_{r}^{k+1}(x)\partial_{x}^{r}u^{k+1}(x)}=f^{k+1}\left(x\right),\label{eq:time-discrete}
\end{eqnarray}
where $f^{k+1}=\mu_{\tau}^{\alpha}\left(u^{k}-\sum_{j=0}^{k-1}a_{k,j}^{\alpha}\left(u^{j+1}-u^{j}\right)\right)+S^{k+1}$. 

So at each time step, we need to solve the higher-order differential
equation (\ref{eq:time-discrete}). 

We consider the following Bernstein Petrov-Galerkin formulation for
(\ref{eq:time-discrete}):

Find $u_{N}\in V_{N}^{0,n}$ such that 
\begin{eqnarray}
 &  & \mu_{\tau}^{\alpha}(u_{N},v_{N})-\sum_{r=0}^{n}{(b_{r}^{k+1}\partial_{x}^{r}u_{N},v_{N})}=(f^{k+1},v_{N}),\quad\forall v_{N}\in W_{N}^{0,n}.\label{eq:Bern-Spectral}
\end{eqnarray}
Using repeated integration by parts along with conditions (\ref{eq:ConditionsOfTest})-(\ref{eq:ConditionsOfTrial}),
we can rewrite (\ref{eq:Bern-Spectral}) as 
\begin{eqnarray}
 &  & \mu_{\tau}^{\alpha}(u_{N},v_{N})-\sum_{r=0}^{n}{(-1)^{\left\lfloor \frac{r+1}{2}\right\rfloor }(\partial_{x}^{\left\lfloor \frac{r}{2}\right\rfloor }u_{N},\partial_{x}^{\left\lfloor \frac{r+1}{2}\right\rfloor }(b_{r}^{k+1}v_{N}))}=(f^{k+1},v_{N}).\label{eq:Bern-spectral-Petrov}
\end{eqnarray}
We expand the approximate solution of (\ref{eq:time-discrete}) in
terms of the basis functions (\ref{eq:trialFunctions}) of $V_{N}^{0,n}$,
i.e., 
\begin{eqnarray}
u^{k+1}(x) & = & \sum_{j=\left\lfloor \frac{n}{2}\right\rfloor }^{N-\lfloor\frac{n+1}{2}\rfloor}{c_{j}^{k+1}\phi_{j}(x)}.\label{eq:approxSolu}
\end{eqnarray}
Choosing the modal functions introduced in Theorem \ref{TheoremModal}
as the test functions, the (\ref{eq:Bern-spectral-Petrov}) is written
equivalently as 
\begin{eqnarray}
 &  & \mathbf{A}\mathbf{c}^{k+1}=\mathbf{f}^{k+1},\label{eq:Systems of equations}
\end{eqnarray}
where $\mathbf{f}^{k+1}=[f_{i}^{k+1}:0\leq i\leq N-n]$ and the $\mathbf{A}$
is given by 
\begin{eqnarray}
 &  & \mathbf{A}=\mu_{\tau}^{\alpha}\mathbf{Q}-\sum_{r=0}^{n}{(-1)^{\left\lfloor \frac{r+1}{2}\right\rfloor }\mathbf{R}_{r}},\label{eq:matrices}
\end{eqnarray}
where 
\begin{eqnarray}
 &  & \mathbf{Q}=[(\phi_{j},\psi_{i})],\,\mathbf{R}_{r}=[(\partial_{x}^{\left\lfloor \frac{r}{2}\right\rfloor }\phi_{j},\partial_{x}^{\left\lfloor \frac{r+1}{2}\right\rfloor }(b_{r}^{k+1}\psi_{i}))],\quad0\leq r\le n,\label{eq:matrices-1}
\end{eqnarray}
for $0\leq i\leq N-n,\,\left\lfloor \frac{n}{2}\right\rfloor \leq j\leq N-\lfloor\frac{n+1}{2}\rfloor.$
The matrices are $(N-n+1)\times(N-n+1)$. The integrals of $\mathbf{f}^{k+1}$
may be approximated by a Gauss-quadrature rule. 

Note that $u^{0}=g$ is given by the initial condition and $u^{k+1},\,k\geq0,$
are obtained from (\ref{eq:approxSolu}) through solving (\ref{eq:Systems of equations})\@.
By the three-term relation (\ref{eq:ThreeTermDer}), (\ref{eq:der of DBPs})
and the biorthogonality system (\ref{eq:Biorthog}), it is found that
the matrices in (\ref{eq:matrices-1}) (so the coefficient matrix$\mathbf{A}$)
are banded for the problems with constant coefficients.

\section{\label{sec:NumericalEx}Numerical experiments}

Here, we provide some numerical examples to illustrate the accuracy
and efficiency of the proposed method.

The errors are measured using the discrete $L^{\infty}$ as 
\begin{eqnarray*}
 &  & L^{\infty}:=\max_{x\in\Omega}{|u(x,T)-u_{N}^{M}(x)|}\approx\max_{0\leq i,j\leq\mathcal{N}}{|u(x_{i},T)-u_{N}^{M}(x_{i})|},
\end{eqnarray*}
where $u$ is the exact solution of the problem, $u_{N}^{M}$ is the
approximation solution at $T=t_{M}=1,$ $x_{i}=\frac{i}{\mathcal{N}}$
and $\mathcal{N}=20$. 

\noindent \textbf{Example 1.} \citep{Gao} Consider the following
time fractional advection-dispersion equation 
\begin{eqnarray*}
 &  & {\partial{}_{t}^{\alpha}}u(x,t)=\kappa_{1}\partial_{x}^{2}u(x,t)-\kappa_{2}\partial_{x}u(x,t)+s(x,t),\quad x\in(0,1),\\
 &  & u(x,0)=g(x),\quad u(0,t)=u(1,t)=0,
\end{eqnarray*}
where $\kappa_{1}$ and $\kappa_{2}$ are the advection and dispersion
coefficients, respectively, and $0<\alpha\leq1$. In this problem
$n=2$, so we choose (\ref{eq:base for n=00003D2}) as the test functions.
The $L^{\infty}$ errors for the method are reported at $t=1$ in
Table \ref{tab:spatialRate} for the case $\kappa_{1}=\kappa_{2}=1$
with exact solution $u=\sin{(2\pi x)\exp{(-t)}}$ and $\tau=0.01$. 

\begin{table}
\begin{centering}
\begin{tabular}{ccccccc}
\hline 
 & \multicolumn{2}{c}{$\alpha=0.25$} & \multicolumn{2}{c}{$\alpha=0.5$} & \multicolumn{2}{c}{$\alpha=0.75$}\tabularnewline
\hline 
$N$ & $L^{\infty}$ & rate & $L^{\infty}$ & rate & $L^{\infty}$ & rate\tabularnewline
\hline 
2 & 4.31E-01 &  & 4.34E-01 &  & 4.37E-01 & \tabularnewline
4 & 5.94E-02 & 2.86 & 5.97E-02 & 2.86 & 6.01E-02 & 2.86\tabularnewline
6 & 3.74E-03 & 6.82 & 3.74E-03 & 6.83 & 3.74E-03 & 6.85\tabularnewline
8 & 1.34E-04 & 11.56 & 1.34E-04 & 11.58 & 1.38E-04 & 11.48\tabularnewline
\hline 
\end{tabular}
\par\end{centering}

\protect\caption{\label{tab:spatialRate}The $L^{\infty}$ error and the spatial rate
of convergence for Example 1.}
\end{table}

\begin{figure}[t]
\begin{centering}
\includegraphics[width=12cm,height=10cm]{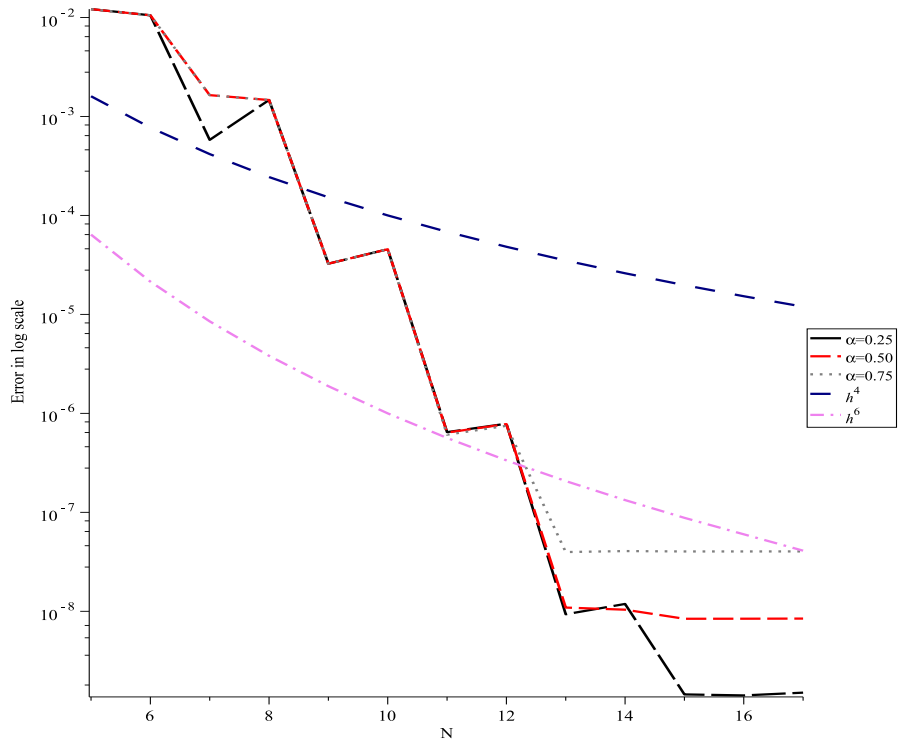}
\par\end{centering}

\protect\caption{\label{fig:Ex2}Convergence of the scheme (\ref{eq:matrices}) and
comparison with the fixed polynomial rates $O(h^{4})$ and $O(h^{6})$.}
\end{figure}

\noindent \textbf{Example 2. }Consider the following equation
\begin{eqnarray*}
 &  & {\partial{}_{t}^{\alpha}}u(x,t)=\partial_{x}u(x,t)+\partial_{x}^{3}u(x,t)-\partial_{x}^{5}u(x,t)+s(x,t),\quad x\in(0,1),\\
 &  & u(x,0)=g(x),\quad u(0,t)=u(1,t)=u_{x}(0,t)=u_{x}(1,t)=u_{xx}(1,t)=0,
\end{eqnarray*}
with the exact solution $u=(1-x)\sin^{2}{(\pi x)}\exp{(-t)}$. Table
\ref{tab:spatialRateEx2} provides the $L^{\infty}$ errors of the
method with $\tau=0.01$ for some fractional orders. The spectral
accuracy of the method is shown in Figure \ref{fig:Ex2} and compared
with fixed rates $O(h^{r}),\,r=4,6$. 

\noindent 
\begin{table}
\begin{centering}
\begin{tabular}{ccccccc}
\hline 
 & \multicolumn{2}{c}{$\alpha=0.25$} & \multicolumn{2}{c}{$\alpha=0.5$} & \multicolumn{2}{c}{$\alpha=0.75$}\tabularnewline
\hline 
$N$ & $L^{\infty}$ & rate & $L^{\infty}$ & rate & $L^{\infty}$ & rate\tabularnewline
\hline 
6 & 1.05E-02 &  & 1.05E-02 &  & 1.05E-02 & \tabularnewline
8 & 1.47E-03 & 6.85 & 1.47E-03 & 6.85 & 1.47E-03 & 6.85\tabularnewline
10 & 4.55E-05 & 15.56 & 4.55E-05 & 15.56 & 4.56E-05 & 15.55\tabularnewline
12 & 7.86E-07 & 22.26 & 7.79E-07 & 22.31 & 7.49E-07 & 22.54\tabularnewline
\hline 
\end{tabular}
\par\end{centering}

\protect\caption{\label{tab:spatialRateEx2}The $L^{\infty}$ error and the spatial
rate of convergence for Example 2.}
\end{table}

\section{\label{sec:Con}Conclusion}

In this paper, utilizing DBPs, a compact formula for the modal basis
functions to solve higher-order BVPs was presented. Using these modal
functions, a Bernstein-spectral Petrov-Galerkin method was established
for a class of time-fractional PDEs. Some numerical examples have
been provided to show the efficiency and spectral accuracy of the
method. The proposed method can be implemented for various time fractional
PDEs on bounded spatial domains.

\end{document}